\documentclass[10pt, a4paper]{amsart}

\usepackage{amsmath,amssymb,amscd,amsfonts}

\newtheorem{thm}{Theorem}[section]
\newtheorem{lem}[thm]{Lemma}
\newtheorem{rem}[thm]{Remark}
\newtheorem{prop}[thm]{Proposition}
\newtheorem{cor}[thm]{Corollary}

\newtheorem{ex}[thm]{Example}

\newcommand\Pic{{\text{\rm Pic}}}
\newcommand\alb{\text{\rm a}}
\newcommand\lra{\longrightarrow}
\newcommand\Alb{{ A}}
\newcommand\Mod{\text{\rm Mod}}
\newcommand\ot{{\otimes}}
\newcommand\OO{{\mathcal{O}}}
\newcommand\Hom{{\mathcal{H}om}}
\newcommand\QQQ{{\mathcal{Q}}}
\newcommand\PP{{\mathcal{P}}}
\newcommand\LL{{\mathcal{L}}}
\newcommand\FF{{\mathcal{F}}}
\newcommand\GG{{\mathcal{G}}}

\newcommand\Q{{\mathbb{Q}}}

\newcommand\CC{{\mathbb{C}}}
 
 \newcommand\II{{\mathcal{I}}}
 \newcommand\sh{{\hat {\mathcal S}}}
 \newcommand\s{{ {\mathcal S}}}
 \title{A derived category approach to generic vanishing}
 \author {Christopher D. Hacon}
 \date{}
 \begin{document}
 \abstract{We prove a Generic Vanishing Theorem 
 for coherent sheaves on an abelian variety over an algebraically closed
 field $k$. When $k=\CC$ this implies a conjecture of Green and Lazarsfeld.}
 \endabstract
 \maketitle
 \section{\label{intro}Introduction}
 In \cite{GL1} and {\cite{GL2}, Green and Lazarsfeld prove the following result which is
an essential tool in the study of irregular varieties:
 \begin{thm} (Generic Vanishing Theorem.) Let $X$ be a smooth complex
projective variety. Then
 then every irreducible component of
$$V^i(\omega _X):=\{ P\in \Pic ^0(X)\ |\ h^i(X,\omega _X \ot P)\ne 0 \}$$
is a translate of a subtorous of $\Pic ^0(X)$ of codimension at least
$$i-\left(\dim (X)-\dim (\alb _X(X))\right).$$ 
If $\dim (X)=\dim \alb _X(X)$ then
there are inclusions:
$$V^0(\omega _X )\supset V^1(\omega _X )\supset ...
\supset V^{\dim (X)}(\omega _X)=
\{\OO _X\}.$$ 
\end{thm}
The proof of this theorem (and of its generalizations) 
relies heavily on the use of
Hodge theory. It is a natural question to try and understand 
to what extent these
results depend on Hodge theory, and what aspects of the proof can be replaced
by an algebraic approach (cf. \cite{EV} 13.13.d). Since the 
above theorem can also
be understood in terms of the sheaves $R^i \alb _{X,*}(\omega _X)$, 
it is also natural to
ask to what extent the Generic Vanishing Theorem generalizes to coherent 
sheaves on an abelian variety. It is surprising that both questions can be 
answered via simple algebraic methods. (Transcendantal methods are of course
used in the proofs of Koll\'ar's theorems on higher direct images of
dualizing sheaves. See however \cite{EV} for a discussion of
algebraic approaches to vanishing theorems.)

 We begin by fixing some notation.
 Let $A$ be an abelian variety over an algebraically
 closed field $k$, $\hat{A}=\Pic ^0(A)$ the dual
 abelian variety and 
 $\LL$ be the normalized Poincar\'e line bundle. 
 We will denote by 
 $p_A,p_{\hat {A}}$ the projections of $A\times \hat{A}$ on to $A,\hat{A}$.
 For any ample line bundle $L$ on $\hat{A}$,
 the isogeny $\phi _L: \hat {A}\longrightarrow A$ is defined by 
$\phi _L(\hat {a})=t_{\hat {a}}^* L^\vee \ot L$. 
Let $\hat{L}$ be the vector bundle on $A$ defined by 
 $$\hat {L}:=p_{A,*}(p_{\hat {A}}^*L\ot \LL).$$
 One has that $$\phi _L^*\left( \hat {L}^\vee \right) \cong \bigoplus _{h^0(L)}L.$$
 We will prove the following:
 \begin{thm} \label{maint}
 Let $\FF$ be a coherent sheaf on an abelian variety $A$.  
 The following are equivalent:
 \begin{enumerate}
 \item 
 For any sufficiently 
 ample line bundle
 $L$ on $\hat{A}$,  
 $$H^i(A, \FF \ot \hat {L}^\vee )=0\ \ \ \forall \ i>0;$$
 \item 
 There is an isomorphism
 $$R p_{\hat {A},*} (p_A^*D_A(\FF ) \ot \LL )\cong R^0 p_{\hat {A},*} (p_A^*
 D_A(\FF )\ot \LL ).$$
 \end{enumerate}
 \end{thm}
 We refer to this result as a Generic Vanishing Theorem as it implies in particular
 that every irreducible component of  $$V^i(A,\FF \ot P):=\{P\in \hat{A}\ |\ h^i(A,\FF \ot P)\ne
0\}$$ has codimension at least $i$ in $A$
 (cf. Corollary \ref{C1}) and that one has inclusions
$$V^0(\FF )\supset V^1(\FF )\supset ...\supset V^n(\FF ).$$ 
 When $X$ is a smooth complex projective variety and $\FF =R^i\alb _{X,*}\omega _X$,
one recovers a generalization of the results of Green and Lazarsfeld.
  It is worthwhile to point out that contrary to what one might expect from the
 results of Green and Lazarsfeld, the support of the sheaves $R^0 p_{\hat {A},*}
 (p_A^*D_A(\FF )\ot \LL)$ are not necessarily reduced subvarieties of $\hat {A}$.
 More precisely, we have:
 \begin{ex} Consider a nontrivial extension
 $$0\lra \OO _A \lra V\lra \OO _A\lra 0.$$
 Let $\FF =V$, then $\FF$ satisfies the hypothesis of Theorem \ref{maint}, 
 and one has that 
 $$R^0 p_{\hat {A},*} (p_A^*D_A(V ) \ot \LL )$$
is an Artinian $\OO _{\hat {A},\hat {0}}$ module of length 2.
 \end{ex}
It is also not the case that the loci
$$V^i(\FF ):=\{ P\in \hat {A}\ s.t.\ H^i(\FF \ot P)\ne 0\}$$
(or equivalently that the sheaves $R^i\sh (\FF ):=R^ip_{\hat {A},*}(p_A^*\FF \ot \LL)$)
should be supported on subtori of
$\hat{A}$ as is illustrated by the following:
\begin{ex} 
Let $(A,\Theta )$ be a principally polarized abelian variety over $\CC$, then 
for any ample line bundle $L$ on $\hat {A}$, one has that 
$\OO _A(\Theta )\ot \hat {L}^\vee$ is globally generated (cf. \cite{PP}) and has
vanishing higher cohomology groups. Fix any point $x\in A$ and let $\FF :=
\OO _A(\Theta )\ot \II _x$. From the short exact sequence
$$0\lra \FF \ot \hat {L}^\vee\lra \OO _A(\Theta )\ot \hat {L}^\vee\lra
\OO _A(\Theta )\ot \hat {L}^\vee\ot \CC (x)\lra 0,$$
one sees that $H^i(A,\FF \ot \hat {L}^\vee )=0$ for all $i>0$.
One has $R^i\sh (\FF)=0$ for all $i\geq 2$
and there is an exact sequence 
$$0\lra R^0\sh (\FF) \lra \widehat{\OO _A(\Theta)}\lra P_x \lra R^1\sh (\FF )\lra 0,$$
where $P_x$ is the topologically trivial line bundle on $\hat {A}$
determined by the point $x\in \hat {\hat {A}}=A$ and $\widehat {\OO _A(\Theta )}
:=R^0\sh(\OO _A(\Theta ))$
is identified with $\OO _A(\Theta )^\vee$ 
under the isomorphism $\phi _\Theta :\hat {A}\lra A$.
Since the homomorphism $\widehat {\OO _A(\Theta )}\lra P_x$ is non zero,
one has that it is an inclusion of line bundles.
Therefore, $R^0\sh (\FF )=0$ and
 the sheaf $R^{1}\sh (\FF )$ is supported on a divisor in $|\Theta \ot P_x|$
and hence not on a subtorous of $\hat {A}$.
\end{ex}

 The statement of Theorem \ref{maint}, appears to be very technical, however
 it seems to have concrete applications to the birational geometry of
 complex projective varieties. 
 In particular, we prove Theorem \ref{conjGL}, (a more general version) of
 a conjecture of Green and Lazarsfeld (cf. Problem 6.2 \cite{GL2}):
 \begin{thm} If $X$ is a smooth complex projective variety with
$\dim (\alb _X(X))=\dim (X)$, then for the universal family of
 topological trivial line bundles $\PP \lra X\times \Pic ^0(X)$,
 one has that  $$R^i{\pi_{\hat {A},*}}(\PP )=0$$
 for all $i<\dim (X)$.
 \end{thm}
It should be noted that the results of Green
and Lazarsfeld hold for $X$ a K\"ahler manifold. The methods
of this paper only apply to the projective setting and so the above theorem
does not completely answer Problem 6.2 of \cite{GL2}. It should however
be possible to answer Problem 6.2 using results of K. Takegoshi \cite{Ta}
by using the same methods of this paper. We do not pursue it in this paper.

 In section \ref{Ueno}, we give an application to the higher direct images of
tensor powers of the canonical line bundle under a surjective
morphism to an abelian variety. To be more precise, let
$a:X\lra A$ be a surjective morphism with connected fibers from
 a smooth complex projective variety with $\kappa (X)=0$ to an abelian
variety. In \cite{ChH2} it was shown that if $P_1(X)\ne 0$, then 
$a_*(\omega _X)=
\OO _A$. Here we extend this result as follows: 
{\it For all $N>0$ there exists
a unipotent vector bundle $V_N$ and an inclusion
$V_N\hookrightarrow a_* (\omega _X ^{\ot N})$ which is generically 
an isomorphism
and induces an isomorphism on global sections.} In particular, the sheaves 
$a_*(\omega _X^{\ot N})$ are semi-positive (cf. \cite{V1}). 
It would be interesting to 
see if one can show that the inclusion is an isomorphism in codimension 1, 
and if
this result can be used to give bounds on $\kappa (F_{X/A})$ 
(compare \cite{Ka2}).

 The techniques used to prove Theorem \ref{maint} are standard results
 in the theory of derived categories for which we refer to \cite{Ha2}.
 It is surprising that one is able to obtain such a generalization
 of the results of Green and Lazarsfeld by means of such a simple
 argument. It is our hope that these methods might also give interesting
 applications to other Fourier-Mukai isomorphisms.
 
\medskip
\noindent {\bf Acknowledgments.} The author would like to thank A. Craw for 
useful suggestions. The author was partially supported by NSA grant no:
MDA904-03-1-0101 and by a grant from the Sloan Foundation.
Shortly after the submission of this paper, the author learned that
G. Pareschi had also obtained a proof of Problem 6.2 of \cite{GL2}
by a completely different approach. 
 \subsection{Notation and conventions}
 We will work over $k$ an algebraically closed field. On a smooth variety,
 we will identify Cartier divisors and line bundles, and we will use the additive
 and multiplicative notation interchangeably. 
 For any normal variety, with structure map $f:X\lra Spec (k)$,
 one has a dualizing complex $\omega ^\centerdot _X:=f^!k$ (cf. \cite{Ha2}).
If $X$ is a smooth projective variety, 
then $K_X$ will be a canonical divisor and $\omega _X=\OO _X(K_X)$, and we denote by
$\kappa(X)$ the Kodaira dimension,  by
$q(X):=h^1(\OO_X)$ the {\em irregularity} and by $P_m(X):=h^0(\omega _X^{\ot m
})$ the {\em $m-$th plurigenus}.
If $f\colon X\to Y$ is a
morphism of smooth projective varieties, we write $K_{X/Y}:=K_X-f^*K_Y$.  
We denote by  $\alb \colon X\to \Alb$ the Albanese map  and by 
$\hat{\Alb }=\Pic ^0(\Alb )$ the dual abelian variety to $\Alb $ which parameterizes all
topologically trivial line bundles on $\Alb$. Recall that $\Pic ^0(\Alb )=\Pic ^0(X)_{red}$.
For a $\Q -$divisor $D$
we let $\lfloor D\rfloor$ be the integral part and
$\{D\}$ the fractional part. Numerical
equivalence is denoted by $\equiv$ and we write $D\prec E$ if
$E-D$ is an effective divisor. 
If $L$ is a Cartier divisor, $|L|$ denotes the complete linear series associated
to $L$. The rest of the notation is
standard in algebraic geometry.
 
 \section{\label{pre} Preliminaries}
 \subsection{Derived Categories}

 For a covariant left exact functor $F:A\lra B$ of abelian categories $A,B$,
 $RF:D(A)\lra D(B)$ will denote the 
 right derived functor between the corresponding 
 derived categories. 
 In particular for any $X\in Ob(A)$, the $i$-th cohomology
 group of $RF (X)$ is just $R^iF(X)$. For any scheme $X$, 
 $D(X)$ is the derived category
 of the category of $\OO _X$-modules $\Mod (X)$.
 $D_c(X),D_{qc}(X)$ denote the full subcategories of $D(X)$ consisting 
 of complexes whose cohomologies are coherent, quasi-coherent.
 $D^-(X),D^+(X),D^b(X)$ denote the full subcategories of $D(X)$ consisting
 of complexes bounded above, below, on both sides.
 $F [n]$ denotes the complex obtained by shifting the complex
 $n$ places to the left.
 For the convenience of the reader we recall the following facts:
 \begin{enumerate}
 \item {\bf Projection formula (P.F.).} (cf. \cite{Ha2} \S II) 
 Let $f:X\lra Y$ be a proper
 morphism of quasi-projective varieties. 
 Then there is a functorial isomorphism:
 $$R f_*(F) \ot _Y G\lra R f_*(F \ot _X
 L f^*G )$$
 for $F \in D^-(X)$ and $G \in D^-_{qc}(Y)$.
 
 \item {\bf Grothendieck Duality (G.D.).} (cf. \cite{Ha2} \S VII)
 For an $n$-dimensional variety $X$, 
 the dualizing functor is defined by $D_X(?)=R \Hom (?,\omega ^\centerdot _X[n])$. 
 In particular $\OO _X \cong D_X (\omega ^\centerdot _X [n])$.
 Let $f:X\lra Y$ a morphism of projective varieties, $F\in D^b_{qc}(X)$
 then $$R f_* D_X(F  )\cong D_Y R f_* (F ).$$
 In particular one has that 
$$R \Gamma (D_X F )\cong D_k (R \Gamma (F ) ),$$ 
 (where $R \Gamma=R \gamma _*$ and $\gamma:X\lra k$ is
the structure morphism).
 \end{enumerate}
 \subsection{Higher direct images of dualizing sheaves}
 Recall the following results
 \begin{thm} \label{kollar} {(\cite {Ko1}, \cite{Ko2})}
 Let $X,Y$ be complex projective varieties of dimension $n,k$
 with $X$ smooth.
 Let $f:X\lra Y$ be a surjective map and $L$ an ample line bundle on $Y$. Then
 \begin{enumerate}
 \item $R^if_* \omega _X $ is torsion free for $i\geq 0$,
 \item $H^j(Y,L\ot R^if_*\omega _X )=0$ for $j>0$,
 \item $R f_* \omega _X \cong \sum R^if_* \omega _X [-i].$
 \item $R f_* \OO _X \cong \sum {}^dR^if_* \OO _X$, where 
 ${}^dR^if_* \OO _X=D_Y(R^{n-k-i}f_* \omega _X [k+i])$.
 \end{enumerate}
 \end{thm}

\begin{thm} \label{omegavan} (\cite{Ko3} \S 10) Let $f\colon X\to  Y$ be a  surjective 
map of projective
varieties, $X$ smooth, $Y$ normal. Let $M$ be a line bundle on $X$
such that $M\equiv N+f^* L +\Delta$, where $N$ is a $\Q$-divisor on $X$ which is either
nef and big or numerically trivial, $L$
is a $\Q-$divisor on $Y$
and $\Delta$ has normal crossing support with 
$\lfloor \Delta \rfloor=0$. Then
\begin{enumerate}
\item $R^jf_*(\omega _X\ot M) $ is torsion free for $j\geq 0$;

\item If $L$ is nef and big, and $g:Y\lra Z$ is any morphism with $Z$ projective,
then $H^i(Y,R^jf_*(\omega _X\ot M))
=0$ for $i>0$, $j\geq 0$ and
$$H^i(Y,R^jf_*(K_X+M))=H^i(Z,g_*(R^jf_*(K_X+M))),\ \ i,j\geq 0.$$
\end{enumerate}
\end{thm}

 \subsection{Sheaves on abelian varieties}
 Let $A$ be a $g$-dimensional abelian variety and $\hat{A}=\Pic ^0(A)$ its 
dual abelian variety. We will denote by $p_A,p_{\hat {A}}$ the projections
of $A\times \hat{A}$ on to $A,\hat {A}$. Let $\LL $ be the normalized Poincar\'e 
line bundle on $A\times \hat {A}$.
The main tool in understanding sheaves on an abelian variety is the 
Theory of Fourier-Mukai transforms. Recall that Mukai defines the functor $\sh $ of $\OO_A$-modules
 into the category of $\OO_{\hat {A}}$-modules by
 $$\sh (M)=p _{\hat{A},*}(\LL \ot p_A ^*M).$$ Similarly 
 one defines $\s (N)=p _{{A},*}(\LL \ot p_{\hat{A}} ^*N).$ 
 By \cite{Mu}, one has: 
 \begin{thm} \label{mukai} {\bf Mukai} There exists isomorphisms of functors
 $$R \sh \circ R \s\cong (-1_{\hat{A}})^*[-g]$$
 and
 $$R \s \circ R \sh \cong (-1_{{A}})^*[-g].$$
 \end{thm}
 We will need the following results analogous to cohomology and base change
(\cite{Ha1} III.12 and \cite{EGA}
\S 7):
 \begin{lem} Let $F \in D^b_c(A)$, and $P_{\hat{a}}$ be the topologically
 trivial line bundle corresponding to the point $\hat{a}\in \hat {A}$. 
If $R^i\Gamma (F  \ot P_{\hat {a}})=0$, then
$R^i\sh (F )\ot k (\hat{a})=0$ and the natural map
$$\varphi ^{i-1}(\hat {a}):R^{i-1}\sh (F )\ot k (\hat{a})\lra R^{i-1}\Gamma 
 (F \ot P_{\hat{a}})$$
 is surjective.
 \end{lem}

 We recall for future reference the following example of Mukai (cf. \cite {Mu}
 Example 2.9). 
 \begin{ex} \label {exM} {\bf (Mukai)}
 A vector bundle $U$ on $A$ is unipotent if there exists a filtration
 $$0=U_0\subset U_1\subset ...\subset U_{n-1}\subset U_n=U$$
 such that $U_i/U_{i-1}\cong \OO _A$ for all $1\leq i\leq n$. One has that
 $R\sh (U)=R^g\sh (U)$ and $R^g\sh (U)$ is supported on
 $\hat {0}\in \hat{A}$. This gives an equivalence between the category of
 unipotent vector bundles on $A$ and the category of coherent sheaves on
 $\hat {A}$ supported on $\hat{0}$ i.e. the category of Artinian $\OO 
 _{\hat{A},\hat{0}}$-modules.
 \end{ex}
 \section{Main result}
 We now proceed to prove Theorem \ref{maint}:
 \begin{proof}
 By Grothendieck Duality (G.D.) and the projection formula (P.F.) it follows that
$$ D_k (R \Gamma (\FF \ot 
 \hat {L}^\vee ))\cong^{G.D.}R \Gamma (D_A(\FF \ot \hat{L}^\vee ))\cong$$ 
 $$R \Gamma (D_A(\FF )\ot \hat{L})\cong 
 R \Gamma (D_A(\FF )\ot p_{A,*}(\LL \ot p_{\hat{A}}^*L))\cong ^{P.F.}$$
 $$R \Gamma (L p_A^* D_A (\FF) \ot \LL \ot p_{\hat{A}}^*L)\cong^{P.F.}
R \Gamma (R \sh (D_A(\FF ))\ot L).$$
In particular, $D_k (R \Gamma (\FF \ot 
 \hat {L}^\vee ))$ is a sheaf if and only if $R \Gamma (R \sh
 (D_A (\FF)) \ot L)$ is a sheaf.
The theorem now follows from the following remarks:
\begin{enumerate} 
\item Notice that there is an isomorphism
$$H^0(A,\FF \ot \hat {L}^\vee )^\vee \cong D_k (R \Gamma (\FF \ot 
 \hat {L}^\vee ))$$ 
if and only if $$H^i(A,\FF \ot \hat {L}^\vee )=0\ for\ all\ i>0.$$
 
\item For $L$ sufficiently ample, the coherent sheaves
 $R^j\sh (D_A (\FF))\ot L$ are globally generated and have vanishing higher cohomology. 
Consider the spectral sequence
 $$E_2^{i,j}=R^i\Gamma (R^j\sh (D_A (\FF ))\ot L)\Rightarrow R^{i+j}\Gamma (
R \sh (D_A(\FF )) \ot L).$$
 Since $E_2^{i,j}=0$ for all $i\ne 0$, this sequence degenerates at the $E_2$
 level and hence $$H^0(\hat{A}, R^j\sh (D_A(\FF ))\ot L)\cong 
 R^j\Gamma (R\sh ( D_A(\FF ))\ot L ).$$
Therefore, 
$$R \Gamma (R\sh ( D_A (\FF))\ot L)\cong
R^0 \Gamma (R \sh (D_A (\FF)) \ot L)$$
if and only if 
$$R \sh (D_A(\FF ))\cong R^0\sh (D_A(\FF )).$$
\end{enumerate}
 \end{proof}
\begin{cor} Let $\FF$ be a coherent sheaf as above. For all $i>0$, 
the support of $R^i\sh (\FF )$ 
has codimension at least $i$ in $\hat {A}$.
\end{cor}
\begin{proof} By \cite{Mu} (3.8), one has that
$$R \sh (\FF) \cong 
D_{\hat{A}} \left((-1_{\hat {A}}^*\circ R^0 \sh \circ D_A)(\FF)[g]\right)\cong
R\Hom \left( \GG,\OO _{\hat{A}} \right),$$
where $\GG:=-1^*_{\hat{A}}( R^0\sh ( D_A(\FF )))$ is a sheaf.
It suffices therefore to show that for any coherent sheaf $\GG$ on a smooth 
affine variety $Y=Spec(A)$, one has that
the sheaves $R^i\Hom( \GG, \OO _Y)$ are supported in codimension at 
least $i$.
Let $M=\Gamma (Y,\GG )$ and $W$ be an irreducible component of the support of 
$R^i\Hom ( \GG, \OO _Y)={Ext^i(M,A)}^\backsim$. Let $P\in Spec (A)$ such that
$W=Spec (A/P)$. 
Since localization is exact, one has that
$$0\neq Ext^i (M,A)\ot   A_P\cong
Ext^i (M_P,A_P).$$
Since $A$ is regular, it follows that $A_P$ is regular and hence that
$$i\leq\dim A_P=\dim (Y)-\dim (W)$$ (cf. \cite{Ha1} \S III.6).
\end{proof}
 \begin{cor}\label{C1} Let $\FF$ be a coherent sheaf as above, and $P\in \hat{A}$,
 then:
 \begin{enumerate}
 
 \item If $i\geq 0$ and $H^i(A,\FF \ot P)=0$, then $H^{i+1}(A,\FF \ot P)=0$.

\item For all $i>0$ any irreducible component of the loci
$$V^i(\FF ):=\{ P\in \Pic ^0(A)\ |\ h^i(\FF \ot P)\ne 0\}$$
has codimension at least $i$ in $\hat {A}$.
 
 \item If $H^0(A, \FF \ot P)=0$ then $R^0\sh (D_A(\FF))\ot k(P^\vee )=0$.

 \item If $H^0(A,\FF \ot P)=0$ for all $P\ne \OO _A$, then $\FF$ is a 
 unipotent vector bundle.
 \end{enumerate}
 \end{cor}
 \begin{proof} 
 (1) Since
$$0=H^i(A,\FF \ot P)^\vee \cong H^{-i}(D_\CC R
 \Gamma (\FF \ot P))\cong R^{-i}\Gamma (D_A(\FF )\ot P^\vee ),$$
by Cohomology and Base Change, one sees that the natural
homomorphism
$$0=R^{-i-1}\sh (D_A(\FF ))\ot k(P^\vee)\lra R^{-i-1}\Gamma (D_A(\FF)\ot P^\vee)$$ 
is surjective and hence that $H^{i+1}(A,\FF \ot P)=0$.

(2) Let $W$ be an irreducible component of $V^i(\FF)$ and $P$ a general point in $W$.
 Let $j$ be the largest integer such that $W\subset V^j(\FF )$.
By Cohomology and Base Change (for sheaves), one has that $R^j\sh (\FF )\ot k(P)\ne 0$
and by the previous corollary, one has that $W$ has codimension at least $j$.

 (3) This follows immediately from Cohomology and Base Change.

 (4) By (3), we have that $R^0\sh (D_A(\FF ))$ is supported on 
 the closed point $\hat{0}\in \hat {A}$ corresponding to $\OO _A$.
 By Example \ref{exM}, $R^0 \sh (D_A(\FF ))$ is an Artinian
 module and hence
 $$D_A(\FF )\cong (-1_A)^*R \s R \sh (D_A (\FF ))$$
 is (the shift of) a unipotent vector bundle say $V$. Finally it is easy to see
 that $$\FF=D_A(D_A(\FF))=V^\vee$$ is also a unipotent vector bundle.
 \end{proof}
\section{\label{GL}A conjecture of Green and Lazarsfeld}
When $X$ is a smooth complex projective variety and $\alb :X\lra A$ is the
 Albanese morphism, one recovers a generalization of the conjecture of
 Green and Lazarsfeld mentioned in the introduction. 
 To fix the notation, let $\PP=(\alb ,id_{\hat{A}})^*\LL$ and let $\pi _{\hat A}$ denote
 the projection of $X\times \hat {A}$ onto $\hat{A}$.
 \begin{thm} \label{conjGL}
 Let $X$ be a smooth complex projective variety of dimension 
 $n$, and Albanese dimension $k:=\dim \alb (X)$. 
 Then $$R \pi _{\hat {A },*}(\PP )\cong \sum _{i=0}^{n-k}R^{k+i}
 \sh ({}^dR^i\alb _* \OO _X).$$
 In particular, when $X$ is of maximal Albanese dimension (i.e. $k=n$),
 one has that $R \pi _{\hat {A},*}(\PP )\cong R^n\sh ({}^d R^0\alb_* \OO _X)$ is a 
 sheaf.
 \end{thm}
 \begin{proof} 
 Let $Y:=a(X)$. We have that
 $$R \pi _{\hat {A},*}(\PP )\cong R \pi _{\hat {A},*}((\alb \times id_{\hat{A}})^*
 \LL)\cong ^{P.F.}
 R p_{\hat {A},*}( L p_A^* (R \alb _* \OO _X)\ot \LL )\cong$$
 $$R p_{\hat {A},*}( L p_A^*
 (\sum _{i=0}^{n-k}{}^dR^i\alb _* \OO _X)\ot \LL)
 \cong 
 \sum _{i=0}^{n-k}R \sh ({}^dR^i\alb_* \OO _X).$$
 Since $${}^dR^ia_* \OO _X\cong D_Y(R^{n-k-i}\alb _* \omega _X[k+i])\cong
 D_A(R^{n-k-i}\alb _* \omega _X[k+i])$$
 it suffices to show that the sheaves $\FF ^{n-k-i}:= R^{n-k-i}\alb _* \omega _X$
 satisfy the hypothesis of Theorem \ref{maint}. 
 This is however an immediate consequence of Theorem \ref{kollar} and the  
observations that follow:

 Let $\FF =\FF ^{n-k-i}$. As mentioned in the Introduction,
 $\hat {L}:=R^0\s (L)\cong R \s (L)$ is a vector bundle on $A$ of rank
 $h^0(L)$ and
 $$\phi _L^* \left( \hat{L}^\vee \right) \cong \left( L \right) ^{\oplus h^0(L)}.$$
 We have that 
 $H^i(A,\FF \ot \hat {L}^\vee )$ is a direct summand of 
 $$H^i(A,\phi _{L,*}(\OO _{\hat {A}})\ot \FF \ot \hat {L}^\vee  )\cong ^{P.F.}H^i(\hat {A},
 \phi _L^* (\FF \ot \hat {L}^\vee ))\cong \bigoplus _{i=1}^{h^0(L)}H^i(\hat {A},
 \phi _L^* (\FF )\ot L ).$$ So it suffices to show that $H^i(\hat {A},
 \phi _L^* (\FF )\ot L )=0$ for all $i>0$.
 Let $X^\prime:=X\times _A \hat {A}$ and $\alb ^\prime:X^\prime
\lra \hat {A}$ be the induced morphism.
 By flat base change, $\phi _L^* (\FF )\cong R^{n-k-i}{\alb }^\prime_* \omega 
_{X^\prime}$, and so by
 Theorem \ref{kollar}, we have the required vanishing: $$h^i(\hat {A}, 
 R^{n-k-i}{\alb }^\prime_* \omega _{X}^\prime\ot L)=0\ for \ all\ i>0.$$

 \end{proof}
 \begin{cor} Let $\alb:X\lra \Alb$ be as in the preceding theorem, then
\begin{enumerate}
\item Every irreducible component of the loci $V^i(R^ja_*\omega _X)$ has
codimension at least $i$ in $\hat {A}$.

 \item  For all $i<\dim \alb (X)$ and general $P\in \Pic ^0(X)$,
 one has that $H^i(X,P)=0$.

 \item If $H^j(R^i\alb _*\omega _X\ot P)=0$, then $H^l(R^i\alb _*\omega _X\ot P)=0$
 for all $l>j\geq 0$ and $0\leq i\leq k$.
 \end{enumerate} 
\end{cor}
 \begin{proof} This is analogous to Corollary \ref{C1}.
 \end{proof}
 \begin{rem} Following results of Deligne, Illusie and Raynaud (cf. \cite{EV}
 \S 11),
 one can recover some similar statements for $X$ a proper smooth scheme over a perfect field 
 $k$ admitting a lifting $\tilde {X}$ to $W_2(k)$ and $char (k)>\dim (X)$.
 \end{rem}
\section{\label{Ueno}Higher direct images of pluricanonical line bundles}
Recall the following cf. \cite{Mo}:

\noindent {\bf Conjecture K} (Ueno)
  {\em Let $X$ be a smooth complex projective variety with $\kappa
(X)=0$, and let $\alb :X\lra \Alb $ be its Albanese morphism. Then
\begin{enumerate}
\item $\alb$ is surjective and has connected fibers, i.e. $\alb$ is an 
algebraic fiber space;
\item if $F$ is the general fiber of $\alb$, then $\kappa (F)=0$;
\item there is an \'etale covering $B\lra \Alb$ such that 
$X\times _\Alb B$ is birationally equivalent to $F\times B$ over $B$.
\end{enumerate} }

In \cite{Ka1}, Kawamata proved $(1)$ above.
In this section, we will give some evidence towards $(2)$.
Ideally, one would like to show that if $\kappa (X)=0$ and $P_1(X)=1$,
then for all $N>0$, there is an isomorphism $\OO _A \cong 
\alb _*(\omega _X ^{\ot N})^{\vee \vee}$. 
This would imply (2) above and give convincing evidence towards (3). In \cite{ChH2}
this was established for $N=1$. We will show that for $N\geq 2$ there is 
a unipotent vector bundle $V_N$ and an inclusion $V_N\hookrightarrow
 \alb _*(\omega _X^{\ot N})$ 
which is a generic isomorphism.

We begin by recalling some well known properties of multiplier
ideal sheaves. We refer to \cite{La} for a more complete treatment.

Let $X$ be a smooth projective variety and $|L|$ a (non-empty)
linear series on $X$, let $\nu : \tilde {X}\lra X$
be log resolution of $|L|$ i.e. a proper birational morphism from a smooth projective
variety such that $\nu ^*|L|=|M|+F$ where $|M|$ is base point free and
$F\cup \{Exceptional \ locous\ of \nu \}$ has normal crossings support.
For any rational number $c>0$, define the multiplier ideal sheaf
associated to $c$ and $|L|$ by
$$\II (X,c \cdot |L|):= \nu _* (K_{\tilde {X}/X}-\lfloor c F \rfloor ).$$
When $c=1$ and $D=|L|$ this is denoted by  $\II (D)$.
For every $k\geq 1$, one has that $$\II ({c}\cdot |L|) \subset \II
(\frac {c}{k}\cdot |kL|).$$
Therefore, the family of ideals $$\{ \II \left( \frac {c}{k}\cdot |kL|\right)
\} _{k\geq 0}$$
has a unique maximal element which we denote by $\II ( c \cdot || L||)$.
For all $k \gg 0$ one has
$$\II( c \cdot || L||)= \II (\frac {c}{k}\cdot |kL|).$$

Consider now $f:X\lra Y$  a surjective morphism of projective varieties, with 
connected fibers and $X$ smooth. Fix $H$ an ample line bundle on $Y$, and $L$ a
line bundle on $X$ with $\kappa (L)\geq 0$.
For all $t\geq 0$ let
$$\II _t:=\II \left( \frac{1}{t}\cdot ||tL+f^*H||\right) .$$
For any sufficiently big and divisible integer $k$, one has that
(after replacing $k$ by $kt$)
$$\II _t=\II \left( \frac {1}{k}\cdot |kL+\frac {k}{t}f^*H|\right)$$ 
and similarly for $\II _{t+1}$. Let $\nu :\tilde {X}\lra X$
be a log resolution of $|kL+\frac {k}{t}f^*H|$ and of 
$|kL+\frac {k}{t+1}f^*H|$, ($k$ sufficiently big and divisible)
so that $$\nu ^*|kL+\frac {k}{t}f^*H|=|M|+F,\ \ \ \ 
\nu ^*|kL+\frac {k}{t+1}f^*H|=|M'|+F'$$
with $|M|,|M'|$ base point free and $F,F'$ with simple normal
crossings support.
From the inclusion of linear series $$\nu ^*|kL+\frac {k}{t+1}f^*H|
\times \nu ^* f^* |\frac {k}{t(t+1)}H|\lra
\nu ^*|kL+\frac {k}{t}f^*H|$$
one sees that $F\prec F'$
and therefore
$$\II _{t+1}=\nu _* (K_{\tilde {X}/X}-\lfloor \frac {1}{k} F'\rfloor )
\subset \nu _* (K_{\tilde {X}/X}-\lfloor \frac {1}{k} F\rfloor )=\II _t.$$
\begin{prop} There exists an integer $t_0$ such that for all $t\geq t_0$,
one has $$f_*(\omega _X \ot L \ot \II _t)=f_*(\omega _X \ot L \ot 
\II _{t_0}).$$ 
\end{prop}
\begin{proof}
For all $m\geq 1$ and $t\geq 2$,
$$K_{\tilde {X}} +\nu^*L -\lfloor \frac {1}{k}F \rfloor +m\nu^*f^*H
\equiv K_{\tilde{X}}+\frac {1}{k}{M}+\{ \frac {1}{k} F\} +(m-\frac {1}{t})\nu^*f^*H$$
with $\frac{1}{k}M$ nef and big, $(m-\frac {1}{t})H$ ample, $\{(1/k)F\}$
has normal crossings support and $\lfloor \{(1/k)F\}\rfloor=0$.
By Theorem \ref{omegavan}, one has that
$$H^i(Y, f_* (\omega _X \ot L \ot f^* H^{\ot m}\ot \II _t))=$$
$$
H^i(Y, (f \circ \nu)_*(\omega _{\tilde {X}}\ot \nu ^*L \ot \nu ^*f^* H^{\ot m}
(-\lfloor \frac {1}{k}F \rfloor )))=0.$$
Consider now the exact sequence of coherent sheaves on $Y$
$$0\lra f_*(K_X \ot L\ot \II _{t+1})\ot H^{\ot m}
\lra f_*(K_X \ot L\ot \II _t)\ot H^{\ot m}\lra \QQQ _t\ot H^{\ot m}\lra 0,$$
The coherent sheaves $f_* (\omega _X \ot L \ot
f^* H^{\ot m}\ot \II _t)$, 
$f_* (\omega _X \ot L \ot
f^* H^{\ot m}\ot \II _{t+1})$, have
vanishing higher cohomology groups for $m\geq 1$
and therefore, $\QQQ _t\ot H^{\ot m}$ also has vanishing 
higher cohomology groups for all $m\geq 1$.
We may assume that $H$ is very ample.
By Mumford regularity, $\QQQ _t\ot H^{\ot m}$ 
is globally generated for all $m\geq \dim (Y)+1$.
It follows that if $f_*(K_X \ot L\ot \II _t)$ and 
$f_*(K_X \ot L\ot \II _{t+1})$ are not isomorphic 
(i.e. if $\QQQ _t$ is non zero),
then $H^0(\QQQ _t \ot H^{\ot m})\ne 0$. Therefore,
$$h^0(f_*(K_X \ot L\ot \II _{t+1})\ot H^{\ot m})<
h^0(f_*(K_X \ot L\ot \II _{t})\ot H^{\ot m}).$$
In particular, for fixed $m=\dim (Y)+1$, and $t\geq 2$,
one has that $f_*(K_X \ot L\ot \II _t)\ot H^{\ot m}$ and 
$f_*(K_X \ot L\ot \II _{t+1})\ot H^{\ot m}$ 
are not isomorphic for at most finitely many
values of $t$. The Proposition now follows.
\end{proof}
\begin{cor} \label{van} Let $e:\tilde {Y}\lra Y$ be any \'etale morphism.
For all $t\gg 0$, $i>0$
$H,M$ ample on $Y,\tilde{Y}$, one has 
$$H^i(\tilde{Y}, e^*f_*(\omega _X \ot L \ot \II 
(||L+ \frac {1}{t}f^*H||))\ot M)=0.$$
\end{cor}
\begin{proof} For $t\gg 0$, one has that $M\ot  (1/t)e^*H^\vee$ is ample.
Let $\tilde{X}=X\times _Y\tilde {Y}$ and $\epsilon 
:\tilde{X}\lra X$ the induced \'etale morphism.
Proceeding as above, one sees that for $g:\tilde{X}\lra \tilde{Y}$,
$$H^i(\tilde{Y}, e^*f_*(\omega _X \ot L \ot \II _t)\ot M)=
H^i(\tilde{Y}, g_*\epsilon ^* (\omega _X \ot L \ot \II _t)\ot M)=$$
$$H^i(\tilde{Y}, g_* (\omega _{\tilde{X}} \ot \epsilon ^*L \ot g^* M
\ot \epsilon ^* \II _t)).$$
It follows that these cohomology groups 
vanish for all $t\gg 0$.  
\end{proof}
\begin{cor} \label{cm} If $L=\omega _X ^{\ot (N-1)}$ for some $N\geq 2$, then for all $t\gg 0$
one has that:
\begin{enumerate}
\item $H^0(Y,f_*(\omega _X\ot L\ot \II _t))\cong H^0(Y,f_*(\omega _X ^{\ot N}));$
\item For general $y\in Y$, the rank of $f_*(\omega _X\ot L\ot \II _t)$ at $y$ 
is $P_N(F)$.
\end{enumerate}
\end{cor}
\begin{proof} (1) Let $\nu :\tilde {X}\lra X$ be a log resolution of $|NK_X|$,
$|k(N-1)K_X|$ and $|k(N-1)K_X+(k/t)f^*H|$ (for some $t\gg k\gg 0$ sufficiently
big and divisible). The corresponding linear series have
moving parts $|M_N|,|M_{k(N-1)}|,|M|$
which are base point free and the corresponding fixed components
$F_N,F_{k(N-1)},F$ that have simple normal crossings support.
It is easy to see that
$$\frac{F}{k}\prec \frac {F_{k(N-1)}}{k}\prec \frac {N-1}{N}{F_N}\prec F_N.$$
In particular we have inclusions
$$|N\nu ^* K_{{X}}-F_N|\lra |N\nu ^*K_X-\lfloor \frac{F}{k}\rfloor | \lra |N\nu ^*K_X|
\lra |NK_{\tilde {X}}|.$$
Since $h^0( \nu^*\omega _X^{\ot N}(-F_N))=h^0(\nu ^* \omega _X ^{\ot N})=
h^0(\omega _{\tilde {X}}^{\ot N})$,
these are all isomorphisms. Finally, we have that
$$h^0(f_*(\omega _X \ot L \ot \II _t))=h^0(\omega _{\tilde {X}}\ot \nu ^* \omega _X^{\ot (N-1)}
(-\lfloor \frac{1}{k}F\rfloor ))=h^0(\omega _{\tilde {X}}^{\ot N})=P_N(X).$$

To see the (2) recall that for fixed $k\gg 0$ and for some $t\gg 0$,
one has that
$$\FF:=a_*\nu _*(\omega _{\tilde {X}}\ot
\nu ^* \omega _X^{\ot (N-1)}-\lfloor \frac {1}{k}F\rfloor )$$
where $\nu:\tilde {X}\lra X$ is a log resolution of
$|k(N-1)K_X+(k/t)a^*H|$ and so
$$|k(N-1)K_X+\frac{k}{t}\nu ^* a^*H|=|M|+F$$
with $|M|$ base point free and $F$ with simple normal crossings support. 
We may assume that $\nu$ also induces a log resolution of $|NK_{X_y}|$
i.e. that
$$|NK_{\tilde{X}_y}|=|M'|+F'$$ with $|M'|$ base point free and $F'$ a divisor
with simple normal crossings support.
By a result of Viehweg (c.f. \cite{V1}, \cite{V2}, see also \cite{Ko4} \S 10),
we have that $$\lfloor \frac {F}{k}\rfloor \prec F'.$$
In particular, at a generic point $y\in Y$, one has that
$$\FF \ot \CC (y)\cong H^0(\FF \ot \CC (y))\cong
H^0(X_y,\OO _{X_y}(NK_{\tilde X }-\lfloor \frac {F}{k}\rfloor ))\cong $$
$$H^0(X_y, \OO _{X_y}(NK_{{\tilde X}_y}))\cong \CC ^{P_N(X_y)}.$$
\end{proof}
\begin{prop} \label{T1}
Let $X$ be a smooth complex projective variety with $\kappa (X)=0$
and $P_1(X)=1$. Then for all $N>0$, there exists a unipotent vector bundle $V_N$
of rank $P_N(F)$ and an inclusion $V_N\hookrightarrow \alb _*(\omega _X ^ {\ot N})$.
\end{prop}
\begin{proof}
Let $X$ be a smooth projective variety with $\dim (X)=n$, q(X)=q, and
$\kappa (X)=0$. By
the above mentioned result of Kawamata, the Albanese morphism
$\alb :X \lra \Alb$ is an algebraic fiber space.
For a fixed ample divisor $H$ on $\Alb$ and any $t\gg 0$ let
$$\FF :=\alb _*\left(\omega ^{\ot N}_X\ot \II (||(N-1)K_X+\frac {1}{t}\alb^*H||)\right) .$$
Since, $\kappa (X)=0$, it is easy to see that $h^0(\omega ^{\ot N}_X\ot P)=0$ for all
$P\ne \OO _X$ and $h^0(\omega ^{\ot N}_X)=1$ (cf. \cite{ChH2} Lemma 3.1).
In particular, as $\FF$ is a non-zero subsheaf of $\alb_*(\omega _X^{\ot N})$, one has
$h^0(\FF \ot P)=0$ for all $P\neq \OO_X$ and 
$h^0(\FF) \leq 1$. By Corollary \ref{cm}, it follows that
$$\FF  \cong V_N,$$
where $V_N$ is an unipotent vector bundle of rank $P_N(F)$.
\end{proof}

\bigskip
\bigskip
\begin{minipage}{13cm}
\parbox[t]{9.5cm}{Christopher D. Hacon\\
Department of Mathematics,
University of Utah\\
155 South 1400 East, JWB 233\\
Salt Lake City, UT 84112-0090 - USA\\
email: hacon@math.utah.edu\\
fax: 801-581-4148}
\end{minipage}
\end{document}